\documentclass[letterpaper, 10 pt, conference]{ieeeconf}

\IEEEoverridecommandlockouts

\usepackage{amsthm}	
\usepackage{cite}
\usepackage{amsmath,amssymb,amsfonts}
\usepackage{algorithm}
\usepackage{algpseudocode}
\usepackage{graphicx}
\usepackage{textcomp}
\usepackage{xcolor}
\usepackage{graphicx}
\usepackage{textcomp}
\usepackage{eqnarray} 
\usepackage{mathalfa} 
\usepackage{algorithm}
\usepackage{algpseudocode}
\usepackage{adjustbox}                          
\usepackage{graphicx}
\usepackage[utf8]{inputenc}
\usepackage{psfrag}
\usepackage{mathtools}
\usepackage{psfrag}
\usepackage{url}
\usepackage{tabularx}

\usepackage{algorithm,algorithmicx}
\usepackage{algpseudocode}

\newcommand{\sC}{\mathcal{C}}
\newcommand{\sD}{\mathcal{D}}

\newcommand{\sL}{\mathcal{L}}

\newcommand{\sT}{\mathcal{T}}

\newcommand{\R}{\mathbb{R}}
\newcommand{\E}{\mathbb{E}}
\newcommand{\Z}{\mathbb{Z}}

\newcommand{\bv}[1]{\mathbf{#1}}

\newcommand{\plant}[1]{p_{#1}^{(x)}}
\newcommand{\control}[1]{p_{#1}^{(u)}}

\newcommand{\optimalcontrol}[1]{{p_{#1}^{(u)}}^\ast}

\newcommand{\modifiedcontrol}[1]{\bar{q}_{#1}^{(u)}}
\newcommand{\data}{\boldsymbol{\Delta}}

\newcommand{\idealplant}[1]{q_{#1}^{(x)}}
\newcommand{\idealcontrol}[1]{q_{#1}^{(u)}}

\newcommand{\KL}{\text{KL}}

\newcommand{\DKL}{\mathcal{D}_{\KL}}

\newtheorem{Remark}{Remark}
\newtheorem{Theorem}{Theorem}

\newtheorem{Lemma}{Lemma}
\newtheorem{Corollary}{Corollary}
\newtheorem{Assumption}{Assumption}

\newtheorem{Problem}{Problem}

\def\BibTeX{{\rm B\kern-.05em{\sc i\kern-.025em b}\kern-.08em
    T\kern-.1667em\lower.7ex\hbox{E}\kern-.125emX}}
\begin{document}

\title{On a Probabilistic Approach for Inverse Data-Driven Optimal Control\\
}

\author{\'Emiland Garrab\'e, Hozefa Jesawada, Carmen Del Vecchio, and Giovanni Russo$^\ast$
\thanks{E. Garrab\'e and H. Jesawada are joint first authors. H. Jesawada, and C. Del Vecchio are with the Department of Engineering, University of Sannio, Benevento, Italy, E. Garrab\'e and G. Russo are with the Department of Information and Electrical Engineering \& Applied Mathematics, University of Salerno, Italy. Emails: {\tt\small \{jesawada,\ c.delvecchio\}@unisannio.it}, {\tt\small \{egarrabe,\ giovarusso\}@unisa.it}.}%
}       

\maketitle

\begin{abstract}
We consider the problem of estimating the possibly non-convex cost of an agent by observing its interactions with a nonlinear, non-stationary and stochastic environment. For this inverse problem, we give a result that allows to estimate the cost by solving a convex optimization problem. To obtain this result we also tackle a forward problem. This leads to the formulation of a finite-horizon optimal control problem for which we show convexity and find the optimal solution. Our approach leverages certain probabilistic descriptions that can be obtained both from data and/or from first-principles. The effectiveness of our results, which are turned in an algorithm, is illustrated via simulations on the problem of estimating the cost of an agent that is stabilizing the unstable equilibrium of a pendulum. 
\end{abstract}


\section{Introduction}

Inferring the intents of an agent by observing its interactions with the environment is crucial to many scientific domains, with applications spanning across e.g., engineering, psychology, economics, management and computer science. Inverse optimal control/reinforcement learning (IOC/IRL) refers to both the problem and the class of methods to infer the cost/reward driving the actions of an agent by observing its inputs/outputs \cite{ABAZAR2020119}. Tackling this problem is relevant to sequential decision-making \cite{GARRABE202281} and can be useful to design data-driven control systems with humans-in-the-loop as well as incentive schemes in sharing economy settings \cite{Sharing_book}.

In this context, a key challenge in IOC/IRL lies in the fact that the underlying optimization  can become ill-posed even when the environment dynamics is linear, deterministic and the cost is convex. Motivated by this, we propose an approach to estimate possibly non-convex costs when the underlying dynamics is nonlinear, non-stationary and stochastic. The approach leverages probabilistic descriptions that can be obtained directly from data and/or from first-principles. Also, the results allow to obtain cost estimates by solving an optimization problem that we show to be convex.
\footnote{Accepted for presentation at the 62nd IEEE Conference on Decision and Control (CDC), 2023.}

\subsubsection*{Related works}  we briefly survey a number of works related to the results and methodological framework of this paper and we refer to \cite{ABAZAR2020119} for a detailed review on inverse problems across learning and control. As remarked in \cite{506395} IRL has its roots in IOC and these methods were originally developed to find control {\em histories} to produce observed output {\em histories}. It was however quickly noticed that even for simple output histories, the resulting control was often infeasible \cite{506395}. More recently, driven by the advances in computational power to process datasets, IRL methods have gained considerable attention. In \cite{ziebart2008maximum} a maximum entropy-based approach is proposed for stationary Markov Decision Processes (MDPs), which is based on a backward/forward pass scheme (see also \cite{Maximum-Entropy-Multi-Agent} for linear multi-agent games). In \cite{Levine_Locally_Optimal} a local approximation of the reward is used and in \cite{NIPS2011_c51ce410} Gaussian processes are exploited, leading to a method that requires matrix inversion operations and to optimization problems that are not convex in general. Instead, in \cite{6630743} manipulation tasks are considered and path integrals are used  to learn the cost, while in \cite{Finn_Deep_Learning} learning is achieved via deep networks for stationary MDPs.  A model-based IRL approach for deterministic systems is presented in \cite{SELF2022110242} for online cost estimation and \cite{LIAN2022110524} tackles the IRL problem in the context of deterministic multiplayer non-cooperative games. The framework of  linearly solvable MDPs is instead leveraged in \cite{IRL-Todorov} and, while it has the advantages of avoiding solving forward MDPs in each iteration of the optimization and of yielding a convex optimizatiopn problem, it also assumes that the agent can specify directly the state transition. We also recall \cite{Risk-Sensitive}, where a risk-sensitive IRL method is proposed for stationary MDPs assuming that the expert policy belongs to the exponential distribution. Also, in the context of IOC, \cite{NAKANO2023110831} considers stochastic dynamics and proposes an approach to learn the parameter of a control regularizer. The IOC problem for known nonlinear deterministic systems with quadratic cost function in the input is also considered in \cite{Rodrigues}. Finally, as we shall see, in order to obtain our results on the inverse problem we also solve a forward problem that involves optimizing, over probability functions, costs that contain a Kullback-Leibler divergence term. We refer to e.g., \cite{GARRABE202281} for a survey on this class of problems in the context of sequential decision-making across learning and control.

\subsubsection*{Contributions} we introduce a number of results to estimate the possibly non-convex and non-stationary cost of an agent by observing its interactions with the environment, which can be nonlinear, non-stationary and stochastic, and for which just a probabilistic description is known. This probabilistic description can be obtained directly from data. Specifically, by leveraging a probabilistic framework, we give a result that enables to estimate the cost by solving an optimization problem that is convex even when the dynamics is nonlinear, non-stationary and stochastic. In order to obtain our result on the inverse problem, which leverages maximum likelihood arguments, we also tackle a forward problem. This leads to formulation of a finite-horizon optimal control problem with randomized policies as decision variables. For this problem, we find the optimal solution and show that this is a probability mass function with an exponential twisted kernel (this is a class of  policies that is often assumed in works on IRL). Also, we turn our result on cost estimation in an algorithm and its effectiveness is illustrated via simulations on the problem of estimating the cost of an agent that is stabilizing the unstable equilibrium of a pendulum. 

While our results are inspired by works on IRL/IOC, this paper offers a number of key technical  novelties. First, we do not require that the agent can specify its state transitions and we do not assume that the expert policy is stationary. Despite this, our approach leads to an  optimization problem to estimate the cost that we prove to be convex.
Moreover, our approach does not require running and solving forward problems in each iteration of the optimization and it does not require the underlying dynamics to be deterministic.


\section{Mathematical Preliminaries and Problem Formulation}

Sets are in {\em calligraphic} and vectors in {\bf bold}. A random variable is denoted by $\mathbf{V}$ and its realization is $\mathbf{v}$. We denote the \textit{probability mass function}, or simply \textit{pmf}, of $\mathbf{V}$  by $p(\mathbf{v})$ and we let $\sD$ be the convex subset of pmfs. Whenever we take the sums involving pmfs we always assume that the sum exists. The  expectation of a function $\mathbf{h}(\cdot)$ of $\mathbf{V}$ is $\E_{{p}}[\mathbf{h}(\mathbf{V})]:=\sum_{\bv{v}}\mathbf{h}(\mathbf{v})p(\mathbf{v})$, where the sum is over the support of $p(\mathbf{v})$; whenever it is clear from the context, we  omit the subscript in the sum. The {joint} pmf of $\mathbf{V}_1$ and $\mathbf{V}_2$ is denoted by  $p(\mathbf{v}_1,\mathbf{v}_2)$ and the {conditional} pmf of $\mathbf{V}_1$ with respect to $\mathbf{V}_2$ is $p\left( \mathbf{v}_1\mid   \mathbf{v}_2 \right)$. Countable sets are denoted by $\lbrace w_k \rbrace_{k_1:k_n}$, where $w_k$ is the generic set element, $k_1$ ($k_n$) is the index of the first (last) element and  $k_1:k_n$ is the set of consecutive integers between (including) $k_1$ and $k_n$. A pmf of the form $p(\bv{v}_0,\ldots,\bv{v}_N)$ is compactly written as $p_{0:N}$ (by definition $p_{k:k} := p_k(\bv{v}_k)$). We use the shorthand notation $p_{k\mid k-1}$ to denote $p_k(\bv{v}_k\mid\bf{v}_{k-1})$. Also,  functionals are denoted by capital calligraphic characters with arguments within curly brackets. We make use of the Kullback-Leibler (KL \cite{KL_51}) divergence, a measure of the proximity of the pair of pmfs  $p(\mathbf{v})$ and $q(\mathbf{v})$, defined as
$\mathcal{D}_{\KL}\left(p \mid\mid q \right):= \sum_{\bv{v}} p(\bv{v}) \ln\left( {p(\bv{v})}/{q(\bv{v})}\right)$. We also recall here the chain rule for the KL-divergence:

\begin{Lemma}\label{lem:splitting_property}	
Let $\bv{V}$ and $\bv{Z}$ be two (possibly, vector) random variables  and let $f(\bv{v},\bv{z})$ and $g(\bv{v},\bv{z})$ be two joint pmfs. Then, the following identity holds:
\begin{equation}
\begin{split}
    \mathcal{D}_{\KL}
\left( f(\bv{v},\bv{z}) \mid\mid g(\bv{v},\bv{z})  \right) =
\mathcal{D}_{\KL}
\left( f(\bv{v}) \mid\mid g(\bv{v})  \right) + \\ \mathbb{E}_{f(\bv{v})}
\left[	
\mathcal{D}_{\KL} 
\left(f(\bv{z}\mid\bv{y})\mid\mid g(\bv{z}\mid\bv{y}) 
\right)
\right].
\end{split}
\end{equation}
\end{Lemma}

\subsection{Set-up of the Control Problem}\label{sec:main_problem}

We let $\bv{X}_k\in\mathcal{X}\subseteq\Z^n$ be the system state at time step $k$ and $\bv{U}_k\in\mathcal{U}\subseteq\Z^p$ be the control input at time step $k$.  Throughout the paper, the time indexing is chosen so that the control input $\bv{u}_k$ is determined based on information available up to $k-1$ and when input $\bv{u}_k$ is applied, the system transitions from state $\bv{x}_{k-1}$ to state $\bv{x}_k$. We let: (i) $\data_k:=(\bv{x}_{k-1},\bv{u}_k)$ be the input-state data pair collected from the system when it is in state $\bv{x}_{k-1}$ and $\bv{u}_k$ is applied; (ii) $\data_{0:N}:=(\{\data_k\}_{1:N},\bv{x}_N)$ be the dataset over the time horizon $\sT:={0:N}$.  We also denote by $p_{0:N}:=p(\data_{0:N})$  the joint pmf of the dataset. We use the wording {\em dataset} to denote a sequence of input-state data. Sometimes, in applications one has available a collection of datasets, which we term as {\em database} in what follows.

\begin{Remark}
As noted in \cite{Peterka_V_Bayesian_Approach_to_sys_ident_1981}, $p_{0:N}$ is a black box type model that can be obtained directly from the data and does not require assumptions on the underlying dynamics. 
\end{Remark}
We now make the standard assumption that the Markov property holds. Then, $p_{0:N}$ can be conveniently partitioned:
\begin{equation}\label{eqn:behavior}
\begin{split}
p_{0:N} & = p_{0}\left(\bv{x}_0\right)\prod_{k=1}^N p_{k\mid k-1} = p_{0}\left(\bv{x}_0\right)\prod_{k=1}^N 
\plant{k\mid k-1}
\control{k\mid k-1},
\end{split}
\end{equation}
where we used the shorthand notation $\plant{k\mid k-1} := p^{(x)}_k	\left(	\mathbf{x}_k \mid \mathbf{u}_k, \mathbf{x}_{k-1} \right)$, $\control{k\mid k-1}:=p^{(u)}_k	\left( \mathbf{u}_k \mid \mathbf{x}_{k-1} \right)$ and $p_{k\mid k-1}:=\plant{k\mid k-1}\control{k\mid k-1}= p\left(\bv{x}_k,\bv{u}_k\mid\bv{x}_{k-1}\right)$. Also, it is useful to  define the joint pmf $\bar{p}_{k-1:k}:=p_{k-1}(\bv{x}_{k-1},\bv{u}_k)$. We say that \eqref{eqn:behavior} is the probabilistic description of the system.

\begin{Remark}
In (\ref{eqn:behavior}), the pmf $p^{(x)}_k	\left(	\mathbf{x}_k \mid \mathbf{u}_k, \mathbf{x}_{k-1} \right)$ describes in probabilistic terms  the evolution of system, while $p^{(u)}_k	\left( \mathbf{u}_k \mid \mathbf{x}_{k-1} \right)$ is the randomized policy from which,  at time step $k$,  the control input is sampled.
\end{Remark}

We let $c_k:\mathcal{X}\rightarrow\R$ be the cost, at time-step $k$, associated to a given state, $\bv{x}_k$. Then, the expected cost incurred when the system is in state $\bv{x}_{k-1}$ and input $\bv{u}_k$ is applied is given by $\E_{p^{(x)}_{k\mid k-1}}\left[c_k(\bv{X}_k)\right]$.
To address the inverse problem for cost estimation (see Section \ref{sec:cost_estimation_formalization} for the problem statement) we first tackle the following {forward} problem:

\begin{Problem}\label{prob:main}
Given a joint pmf 
\begin{equation}\label{eqn:behavior_target}
\begin{split}
q_{0:N}& := q_{0}\left(\bv{x}_0\right)\prod_{k=1}^N 
q^{(x)}_k	\left(	\mathbf{x}_k \mid \mathbf{u}_k, \mathbf{x}_{k-1} \right)
q^{(u)}_k	\left( \mathbf{u}_k \mid \mathbf{x}_{k-1} \right).
\end{split}
\end{equation}
Find the sequence of pmfs, $\left\{\optimalcontrol{k\mid k-1}\right\}_{1:N}$, such that:
\begin{equation}\label{eqn:main_problem}
    \begin{aligned}
{\left\{\optimalcontrol{k\mid k-1}\right\}_{1:N}}\in & \underset{\left\{ \control{k\mid k-1}\right\}_{1:N}}{\text{arg min}}
    \bigg\{\DKL\left(p_{0:N}\mid\mid q_{0:N}\right)\\ 
    & \ \ \ \ \ + \sum_{k=1}^N \E_{\bar{p}_{k-1:k}}\left[\E_{p^{(x)}_{k\mid k-1}}\left[c_k(\bv{X}_k)\right]\right]\bigg\} \\
   &  s.t.  \  \control{k\mid k-1}\in\sD \ \ \forall k\in\sT.
    \end{aligned}
\end{equation}
\end{Problem}
Throughout the paper we make the following standard
\begin{Assumption}\label{asn:bounded_cost}
The optimal cost of Problem \ref{prob:main} is bounded.
\end{Assumption}
For our derivations, it is also useful to introduce  $\idealplant{k\mid k-1} := q^{(x)}_k	\left(	\mathbf{x}_k \mid \mathbf{u}_k, \mathbf{x}_{k-1} \right)$, $\idealcontrol{k\mid k-1}:=q^{(u)}_k	\left( \mathbf{u}_k \mid \mathbf{x}_{k-1} \right)$ and $q_{k\mid k-1}:=\idealplant{k\mid k-1}\idealcontrol{k\mid k-1} = q\left(\bv{x}_k,\bv{u}_k\mid\bv{x}_{k-1}\right)$. 

As we shall see, the solution of Problem \ref{prob:main} is a sequence of randomized policies.  At each $k$, the control input applied to the system, i.e. $\bv{u}^\ast_k$,  is sampled from  $\optimalcontrol{k\mid k-1}$.
In the cost functional of Problem \ref{prob:main}, minimizing the second term minimizes the expected agent cost, while minimizing the first term amounts at minimizing the discrepancy between $p_{0:N}$ and $q_{0:N}$. Hence, the first term in the cost functional can be thought of as a regularizer,  biasing the behavior of the closed loop system towards the pmf $q_{0:N}$.  
Typically, $q_{0:N}$ acts as passive dynamics \cite{todorov_linearly_mdps,9029512} or expresses desired behavior from demonstration databases \cite{gagliardi2020probabilistic}  
See also \cite{GARRABE202281} for a survey on sequential decision-making problems that involve minimizing this class of cost functionals.
 



\subsection{The Inverse Control Problem}\label{sec:cost_estimation_formalization}

The inverse control problem we consider consists in estimating both the cost-to-go for the agent, say $\bar{c}_k(\cdot)$, and the agent cost $c_k(\cdot)$ given a set of observed states/inputs sampled from $\plant{k\mid k-1}$ and from the agent policy.  
In what follows, we denote by $\hat{\bv{x}}_k$ and $\hat{\bv{u}}_k$ the observed state and control input at time-step $k$.  We also make the following:
\begin{Assumption}\label{asn:inverse_problem}
There exist some $\bv{w}_k:=[w_{k,1},\ldots,w_{k,F}]^T$ such that $\bar{c}_k\left({\bv{x}}_k\right) = -{\bv{w}_k}^T\bv{h}\left({\bv{x}}_k \right)$, where $\bv{h}(\bv{x}_k):=[{h}_1(\bv{x}_k),\ldots,{h}_F(\bv{x}_k)]^T$ and ${h}_i:\mathcal{X}\rightarrow\R$  are known functions, $i=1,\ldots,F$.
\end{Assumption}
In what follows, we say that $\bv{h}(\bv{x}_k)$ is the features vector. The assumption, which is rather common in the literature see e.g., \cite{ziebart2008maximum,Maximum-Entropy-Multi-Agent,6630743,SELF2022110242,IRL-Todorov}, formalizes the fact that the cost-to-go can be expressed as a linear combination of given, possibly nonlinear, {\em features} \cite{Goodfellow-et-al-2016}. With our results in Section \ref{sec:learning_cost} we propose a maximum likelihood estimator for the cost (see e.g., \cite{9661376} for a maximum likelihood framework for linear systems in the context of data-driven control).





\section{Main Results}\label{sec:theory_calculations}

\subsection{Computing the Optimal Policy for Problem \ref{prob:main}}\label{sec:optimal_policy}
With the next result we give the solution to  Problem \ref{prob:main}.
\begin{Theorem}\label{thm:prob_main}
Consider Problem \ref{prob:main} and let Assumption \ref{asn:bounded_cost} hold. Then:
\begin{itemize}
\item[(i)] the problem has the unique solution $\{{\optimalcontrol{k\mid k-1}}\}_{1:N}$ , with 
\begin{equation}\label{eqn:optimal_solution_statement}
\begin{split}
& \optimalcontrol{k\mid  k-1} =\frac{\bar{p}_{k\mid  k-1}^{(u)}\exp\left(-\E_{\plant{k\mid  k-1}}\left[\bar{c}_{k}(\bv{X}_{k})\right]\right)}{\sum_{\bv{u}_{k}}\bar{p}_{k\mid  k-1}^{(u)}\exp\left(-\E_{\plant{k\mid  k-1}}\left[\bar{c}_{k}(\bv{X}_{k})\right]\right)},
\end{split}
\end{equation}
where 
$$
\bar{p}_{k\mid  k-1}^{(u)}:=\idealcontrol{k\mid  k-1}\exp\left(-\DKL\left(\plant{k\mid  k-1}\mid  \mid  \idealplant{k\mid  k-1}\right)\right),
$$
and where $\bar{c}_{k}:\mathcal{X}\rightarrow\R$ is obtained via the  backward recursion
\begin{equation}\label{eqn:backward_recursion}
\begin{split}
& \bar{c}_{k}(\bv{x}_{k})= c_{k}(\bv{x}_{k}) -\hat{c}_{k}(\bv{x}_{k}),\\
&\hat{c}_{k}(\bv{x}_{k})= \ln\left(\E_{\idealcontrol{k+1\mid k}}\left[\exp\left(-\DKL\left(\plant{k+1\mid k}\mid\mid\idealplant{k+1\mid k}\right)\right.\right.\right.\\
& \ \ \ \ \ \ \ \ \left.\left. \left.-\E_{\plant{k+1\mid k}}\left[\bar{c}_{k+1}(\bv{X}_{k+1})\right]\right)\right]\right),\\
&\DKL\left(\plant{N+1\mid N}\mid\mid\idealplant{N+1\mid N}\right)\\ 
& \ \ \ \ \ \ \ +  \E_{\plant{N+1\mid N}}\left[\bar{c}_{N+1}(\bv{X}_{N+1})\right] = 0;
\end{split}
\end{equation}
\item[(ii)] the corresponding minimum is given by:
\begin{equation}\label{eqn:cost_statement}
-\sum_{k=1}^N\E_{\bar{p}_{k-1}}\left[\hat{c}_{k-1}(\bv{X}_{k-1})\right],
\end{equation}
where $\bar{p}_{k-1}:=p_{k-1}(\bv{x}_{k-1})$.
\end{itemize}
\end{Theorem}
\noindent {\em Sketch of the proof.}
The full proof is omitted here for brevity and will be presented elsewhere. We give here a sketch of the proof, which is by induction. 

\noindent{\bf Step $1$.} Consider the cost functional in (\ref{eqn:main_problem}). By means of Lemma \ref{lem:splitting_property}, 
Problem \ref{prob:main} can be recast as the sum of the following two sub-problems:
\begin{subequations}
\begin{equation}\label{eqn:split_1_problem_1}
    \begin{aligned}
    \underset{\{\control{k\mid k-1}\}_{1:N-1}}{\text{min}}
    &\bigg\{\DKL\left(p_{0:N-1}\mid\mid q_{0:N-1}\right) \\ 
    &  +\sum_{k=1}^{N-1} \E_{\bar{p}_{k-1:k}}\left[\E_{p^{(x)}_{k\mid k-1}}\left[c_k(\bv{X}_k)\right]\right]\bigg\} \\
   s.t. & \ \control{k\mid k-1}\in\sD \ \ \forall k\in 1:N-1,
    \end{aligned} 
\end{equation}
\text{and}
\begin{equation}\label{eqn:split_1_problem_2}
    \begin{aligned}
    \underset{\control{N\mid N-1}}{\text{min}}
    &\Big\{\E_{\bar{p}_{N-1}}\left[\DKL\left(p_{N\mid N-1}\mid\mid q_{N\mid N-1}\right)\right.\\ 
    & \hspace{1cm}\left.+ \E_{p_{N\mid N-1}}\left[c_N(\bv{X}_N)\right]\right] \Big\}\\
   s.t. & \ \control{N\mid N-1}\in\sD.
    \end{aligned} 
\end{equation}
\end{subequations} 

Hence, the minimum of (\ref{eqn:split_1_problem_2}) is $\E_{\bar{p}_{N-1}}\left[\sC_N\left\{\optimalcontrol{N\mid N-1}\right\}\right]$, with $\sC\left\{\optimalcontrol{N\mid N-1}\right\}$ being the optimal cost obtained by solving
\begin{equation}\label{eqn:split_1_problem_2_revised}
    \begin{aligned}
    \underset{\control{N\mid N-1}}{\text{min}}
    &\DKL\left(p_{N\mid N-1}\mid\mid q_{N\mid N-1}\right) +  \E_{p_{N\mid N-1}}\left[\bar{c}_N(\bv{X}_N)\right] \\
   s.t. & \ \control{N\mid N-1}\in\sD,
    \end{aligned} 
\end{equation}
 where we set $\bar{c}_N(\bv{x}_N):={c}_N(\bv{x}_N) + \hat{c}_N(\bv{x}_N)$, $\hat{c}_N(\bv{x}_N)=0$. This corresponds to the recursion in (\ref{eqn:backward_recursion}) at $k=N$. 

\noindent{\bf Step $2$.} 
The next step is to show that the problem in (\ref{eqn:split_1_problem_2_revised}) can be conveniently written as
\begin{equation}\label{eqn:split_1_problem_2_revised_for_convex}
    \begin{aligned}
    \underset{\control{N\mid  N-1}}{\text{min}}
    &\bigg\{\E_{\control{N\mid  N-1}}\left[\DKL\left(\plant{N\mid  N-1}\mid  \mid  \idealplant{N\mid  N-1}\right)\right.\\ 
    &+\E_{\plant{N\mid  N-1}}\left[\bar{c}_N(\bv{X}_N)\right]\Big]+ {\alpha}(\bv{x}_{N-1})\bigg\}\\
   s.t. & \ \control{N\mid  N-1}\in\sD,
    \end{aligned} 
\end{equation}
where $\alpha(\bv{x}_{N-1}):= \DKL\left(\control{N\mid  N-1}\mid  \mid  \idealcontrol{N\mid  N-1}\right)$. By studying the second variation of the cost functional, it can be shown that this is strictly convex in the decision variable $\control{N\mid  N-1}$. Hence, since the subset $\mathcal{D}$ is convex, the problem in (\ref{eqn:split_1_problem_2_revised_for_convex}) is a convex optimization problem.

\noindent{\bf Step $3$.} We find the solution to the problem in (\ref{eqn:split_1_problem_2_revised}) by using the equivalent formulation given in (\ref{eqn:split_1_problem_2_revised_for_convex}). Since the problem in (\ref{eqn:split_1_problem_2_revised_for_convex}) is convex with a strictly convex cost functional, the unique optimal solution can be found by imposing the stationarity conditions on the Lagrangian, which is given by:
\begin{equation}\label{eqn:probl_lagrangian}
\begin{split}
 \sL(\control{N\mid  N-1},\lambda_N) & = \E_{\control{N\mid  N-1}}\Big[\DKL\left(\plant{N\mid  N-1}\mid  \mid  \idealplant{N\mid  N-1}\right)+\\
 &\hspace{-1cm}\E_{\plant{N\mid  N-1}}\left[\bar{c}_N(\bv{X}_N)\right]\Big] + \DKL\left(\control{N\mid  N-1}\mid  \mid  \idealcontrol{N\mid  N-1}\right) \\
 &+\lambda_N\left(\sum_{\bv{u}_k}\control{N\mid  N-1}-1\right),
\end{split}
\end{equation}
where $\lambda_N$ is the Lagrange multiplier corresponding to the constraint $\control{N\mid  N-1}\in\sD$. Now, by imposing the first order stationarity conditions on $ \sL(\control{N\mid  N-1},\lambda_N)$, it can be shown that the unique optimal solution is given by:
\begin{equation}\label{eqn:optimal_solution_N}
\begin{split}
&\optimalcontrol{N\mid  N-1} = \frac{\bar{p}_{N\mid  N-1}^{(u)}\exp\left(-\E_{\plant{N\mid  N-1}}\left[\bar{c}_N(\bv{X}_N)\right]\right)}{\sum_{\bv{u}_N}\bar{p}_{N\mid  N-1}^{(u)}\exp\left(-\E_{\plant{N\mid  N-1}}\left[\bar{c}_N(\bv{X}_N)\right]\right)}.
\end{split}
\end{equation}
This is  the optimal solution given in (\ref{eqn:optimal_solution_statement}) for $k=N$, with $\bar{c}_N(\bv{x}_N)$ generated via the backward recursion in (\ref{eqn:backward_recursion}). 
Hence, the minimum for the sub-problem in (\ref{eqn:split_1_problem_2}) is
\begin{equation}\label{eqn:cost_N}
\begin{split}
& \E_{\bar{p}_{N-1}}\left[\sC_N\left\{\optimalcontrol{N\mid  N-1}\right\}\right]  = -\E_{\bar{p}_{N-1}}\left[\hat{c}_{N-1}(\bv{X}_{N-1})\right],
\end{split}
\end{equation}
where
\begin{equation*}
\begin{split}
& \hat{c}_{N-1}(\bv{x}_{N-1})\\
&:= \ln\left(\E_{\idealcontrol{N\mid  N-1}}\left[\exp\left(-\DKL\left(\plant{N\mid  N-1}\mid  \mid  \idealplant{N\mid  N-1}\right)\right.\right.\right.\\
& \hspace{2cm}-\E_{\plant{N\mid  N-1}}\left[\bar{c}_N(\bv{X}_{N})\right]\Big)\Big]\bigg).
\end{split}
\end{equation*}
This is the optimal cost for $k=N$ given in (\ref{eqn:cost_statement}). Next, we make use of the minimum found for the sub-problem (\ref{eqn:split_1_problem_2}) to solve the sub-problem corresponding to $k\in 1:N-1$.

\noindent{\bf Step $4$.} 
It can be shown that the problem in (\ref{eqn:split_1_problem_1}) can be again split as the sum of two sub-problems: one sub-problem for the time-steps up to $N-2$ and a sub-problem for $k=N-1$. Moreover, the latter sub-problem is again independent on the former and has the same structure as (\ref{eqn:split_1_problem_2}). 
Then, following the arguments used in Step $3$, we have that the unique optimal solution for the sub-problem at $k=N-1$ is
\begin{equation}\label{eqn:optimal_solution_N-1}
\begin{split}
& \optimalcontrol{N-1\mid  N-2} \\
& =\frac{\bar{p}_{N-1\mid  N-2}^{(u)}\exp\left(-\E_{\plant{N-1\mid  N-2}}\left[\bar{c}_{N-1}(\bv{X}_{N-1})\right]\right)}{\sum_{\bv{u}_{N-1}}\bar{p}_{N-1\mid  N-2}^{(u)}\exp\left(-\E_{\plant{N-1\mid  N-2}}\left[\bar{c}_{N-1}(\bv{X}_{N-1})\right]\right)}.
\end{split}
\end{equation}
Hence, (\ref{eqn:optimal_solution_N-1}) is the optimal solution given in (\ref{eqn:optimal_solution_statement}) for $k=N-1$, with $\bar{c}_{N-1}(\bv{x}_{N-1})$ obtained via the backward recursion in (\ref{eqn:backward_recursion}). 
We can now draw the desired conclusions.

\noindent{\bf Step $5$.} By iterating Step $4$ we find that, at each of the remaining time steps in the window $1:N-2$, Problem \ref{prob:main} can always be split in sub-problems, where the sub-problem corresponding to the last time instant in the window can be solved independently on the others. 
Hence, the optimal solution for the sub-problem is 
\begin{equation}\label{eqn:optimal_solution_general}
\begin{split}
\optimalcontrol{k\mid  k-1} =\frac{\bar{p}_{k\mid  k-1}^{(u)}\exp\left(-\E_{\plant{k\mid  k-1}}\left[\bar{c}_{k}(\bv{X}_{k})\right]\right)}{\sum_{\bv{u}_{k}}\bar{p}_{k\mid  k-1}^{(u)}\exp\left(-\E_{\plant{k\mid  k-1}}\left[\bar{c}_{k}(\bv{X}_{k})\right]\right)}.
\end{split}
\end{equation}
This is the optimal solution given in (\ref{eqn:optimal_solution_statement}) at time $k$, with $\bar{c}_{k}(\bv{x}_{k})$ obtained from the backward recursion in (\ref{eqn:backward_recursion}). Part (i) of the result is then proved. Moreover, the corresponding optimal cost at time $k$ is $-\E_{\bar{p}_{k-1}}\left[\hat{c}_{k-1}(\bv{X}_{k-1})\right]$. Hence, the optimal cost Problem \ref{prob:main} is $-\sum_{k=1}^N\E_{\bar{p}_{k-1}}\left[\hat{c}_{k-1}(\bv{X}_{k-1})\right]$ and this proves part (ii) of the result. \qed

\begin{Remark}
The optimal solution given in Theorem \ref{thm:prob_main} has an exponential twisted kernel. This class of policies is also known in the literature as soft-max/Boltzmann policies and are often assumed in IRL/IOC works, see e.g., \cite{Risk-Sensitive, ziebart2008maximum,6716965}.
\end{Remark}
\subsection{Estimating the Cost}\label{sec:learning_cost}

Next, we show that the cost $\bar{c}_k(\cdot)$ can be estimated by observing a sequence of states sampled from $\plant{k\mid  k-1}$ when control inputs sampled from $\optimalcontrol{k\mid k-1}$ are applied. This can be useful in settings where one has access to observations of  e.g., an {\em expert}. By estimating $\bar{c}_k(\cdot)$ one can also bypass the computation of the cost-to-go via the backward recursion in Theorem \ref{thm:prob_main}. 
With the next result, we propose an estimator for the cost-to-go. The estimator does not require any linearity/stationary assumption and the underlying dynamics can be stochastic and obtained directly from data.
\begin{Theorem}\label{thm:estimator}
Let Assumption \ref{asn:inverse_problem} hold and let $\hat{\Delta}=\{(\hat{\bv{x}}_{0},\hat{\bv{u}}_{1}),\ldots,(\hat{\bv{x}}_{M-1},\hat{\bv{u}}_{M})\}$ be a sequence of data, with $\hat{\bv{x}}_{k} \sim \plant{k\mid  k-1}$, $\hat{\bv{u}}_{k} \sim \optimalcontrol{k\mid k-1}$ and where $\optimalcontrol{k\mid k-1}$ is  from Theorem \ref{thm:prob_main}. 
 Then, the maximum likelihood estimate for  $\bar{c}_k\left({\bv{x}}_k\right)$, say $\bar{c}_k^\ast\left({\bv{x}}_k\right)$, is given by $\bar{c}_k^\ast\left({\bv{x}}_k\right) =   -{{\bv{w}}^\ast_k}^T{\bv{h}(\bv{x}}_k)$, where $\bv{w}_k^\ast$ is obtained by solving the convex optimization problem
\begin{equation}\label{eqn:prob_log_likelihood_statement}
    \begin{aligned}
    &{\bv{w}}^\ast:=\left[\bv{w}^{*T}_1,\ldots,\bv{w}^{*T}_M\right]^T \in \\
    &\underset{{\bv{w}}}{\text{arg min}}
     \Bigg\{\sum_{k=1}^M\bigg({-\E_{p(\bv{x}_k\mid\hat{\bv{x}}_{k-1},{\hat{\bv{u}}}_k)}\left[{\bv{w}}_k^T{\bv{h}({\bv{x}}}_k)\right]} \\
     & \hspace{2cm} +\ln\Big(\sum_{\bv{u}_k}\modifiedcontrol{k\mid  k-1}(\hat{\bv{x}}_{k-1},{\bv{u}}_k) \\ &\hspace{2.5cm}\exp\left(\E_{p(\bv{x}_k\mid\hat{\bv{x}}_{k-1},\bv{u}_k)}\left[{\bv{w}}_k^T{\bv{h}({{\bv{x}}_k}})\right]\right)\Big)\bigg)\Bigg\},
    \end{aligned}
\end{equation}
and where
\begin{equation}\label{eqn:modified_control}
\begin{split}
    \modifiedcontrol{k\mid  k-1}(\hat{\bv{x}}_{k-1},{\bv{u}}_k):=\left(q({\bv{u}}_k\mid\hat{\bv{x}}_{k-1})\right.\\ 
    &\hspace{-5cm}\left.\exp\left(-\DKL\left(p(\bv{x}_k\mid\hat{\bv{x}}_{k-1},{\bv{u}}_k)\mid\mid q(\bv{x}_k\mid\hat{\bv{x}}_{k-1},{\bv{u}}_k)\right)\right)\right),
\end{split}
\end{equation}
\end{Theorem}
{\em Sketch of the proof.}
The result is based on maximum likelihood, leveraging the structure of the policy of Theorem \ref{thm:prob_main}. Convexity follows from the fact that the feasibility domain is convex and the cost function is a linear combination of the log-sum-exp function and of a linear function (in the decision variables). The proof, omitted here for brevity, will be presented elsewhere. \qed

\begin{Remark}
The problem in \eqref{eqn:prob_log_likelihood_statement} is an unconstrained convex optimization problem with a twice differentiable cost. Constraints on the $\bv{w}_k$'s can be added to capture application-specific requirements, such as {\em dwell-time} constraints. 
\end{Remark}

Next, we propose an estimator when the cost, which we simply denote by $c(\cdot)$, is stationary. The result (the proof of which is omitted here for brevity) implies that the cost can be estimated from a {\em greedy} policy obtained via Theorem \ref{thm:prob_main}. Note that, in this case, the decision variable in the resulting optimization is $\bv{w}_{s}\in\mathbf{R}^F$ rather than $\bv{w}\in\mathbf{R}^{F\times M}$.

\begin{Corollary}\label{clry: estimator}
Let Assumption \ref{asn:inverse_problem} hold and consider $\optimalcontrol{k\mid k-1}$ obtained at each $k$ from Theorem \ref{thm:prob_main} with $N=1$. Further, let the cost be stationary. Then, the maximum likelihood estimate for the cost is ${c}^\ast(\bv{x}_{k}) = {-\bv{w}_{s}^\ast}^T{\bv{h}(\bv{x}}_k)$,  where $\bv{w}_{s}^\ast$ is given by:
 ~\begin{equation}\label{eqn:estimation_stationary}
{\begin{aligned}
\bv{w}_{s}^{\ast} \in \underset{\bv{w}_{s}}{\text{arg min}} & \Bigg\{\sum_{k=1}^M\left(
-\E_{p(\bv{x}_k\mid\hat{\bv{x}}_{k-1},{\hat{\bv{u}}}_k)}\left[\bv{w}_{s}^
T{\bv{h}({{\bv{x}}}_k})\right]\right )\\
 & +\sum_{k=1}^M\ln\left(\sum_{\bv{u}_k}\modifiedcontrol{k\mid  k-1}(\hat{\bv{x}}_{k-1},\bv{u}_k)\right.\\
&\hspace{1cm}\exp\left(\E_{p(\bv{x}_k\mid\hat{\bv{x}}_{k-1},{\bv{u}}_k)}\left[\bv{w}_{s}^T{\bv{h}({{\bv{x}}}_k})\right]\right)\bigg)\Bigg\}.
\end{aligned}}
\end{equation}
 with
$\bv{w}_{s} \in \mathbf{R}^{f}$ and $\modifiedcontrol{k\mid  k-1}(\hat{\bv{x}}_{k-1},\bv{u}_k)$ defined in Theorem \ref{thm:estimator}.

\end{Corollary}

Corollary \ref{clry: estimator} implies that, rather conveniently, the cost can be learned from a greedy policy rather than from the optimal policy. The result  can be also turned into an algorithmic procedure with its main steps given in Algorithm \ref{alg:estimator}. 





\begin{algorithm}[ht]
	\caption{Pseudo-code from Corollary \ref{clry: estimator}} \label{alg:estimator}
	\begin{algorithmic}
			\State \textbf{Inputs:}  observed data $\hat{\bv{u}}_1,\ldots,\hat{\bv{u}}_M$ and $\hat{\bv{x}}_0,\ldots,\hat{\bv{x}}_M$, 
            \State \hspace{1.1cm} $f$-dimensional features vector $\bv{h}(\bv{x}_k)$, 
            \State \hspace{1.1cm} $\plant{k \mid k-1}$,  $\idealplant{k \mid k-1}$, $\idealcontrol{k \mid k-1}$ 
		\State \textbf{Output:} $\bar{c}^\ast\left({\bv{x}}_k\right)$
		\For{$ k = 1$  to $M$}
            \State Compute $\modifiedcontrol{k\mid  k-1}(\hat{\bv{x}}_{k-1},{\bv{u}}_k)$ using \eqref{eqn:modified_control}
            \EndFor
        \State Compute $\bv{w}_{s}^\ast$ by solving the problem in \eqref{eqn:estimation_stationary}
	\end{algorithmic}
\end{algorithm}

\section{Application Example}
We illustrate the effectiveness of our results by considering the problem of stabilizing a pendulum on its unstable equilibrium point. Specifically, given a suitable cost, we first used Theorem \ref{thm:prob_main} to compute the optimal policy and then we leveraged Corollary \ref{clry: estimator} to estimate the cost used in the policy. The pendulum dynamics (only used to generate data) is:

\begin{equation}\label{Pendulum_Dynamics}
    {\begin{aligned}
        \theta_{k} &= \theta_{k-1} + {\omega}_{k-1} dt + W_{\theta} \\
    {\omega}_{k} &
    = {\omega}_{k-1} + \left( \frac{g}{l}\sin(\theta_{k-1})+\frac{u_{k}}{ml^{2}}\right)dt + W_{{\omega}}{,}
    \end{aligned}}
\end{equation}


where $\theta_{k}$ is the angular position, $\omega_{k}$ is the angular velocity and $u_{k}$ is the torque applied on the hinged end. The parameter $l$ is the length of the rod, $m$ is the mass of the pendulum, $g$ is the gravity and $dt$ is the discretization step. Also, $W_{\theta}$ and $W_{{\omega}}$ capture Gaussian noise on the state variables. In our experiments we set $W_{\theta}\sim \mathcal{N}(0,0.05)$ and $W_{{\omega}}\sim \mathcal{N}(0,0.1)$. As in \cite{GARRABE202281} we let $\bv{X}_{k}:= [\theta_{k},\omega_{k}]^T$. Also, $\bv{X}_{k}\in\mathcal{X}$ and $u_{k}\in\mathcal{U}$, where the set $\mathcal{X}:=\left[-\pi,\pi\right]\times[-5,5]$ and $\mathcal{U}:=[-2.5,2.5]$ are discretised in $50\times50$ and $20$ bins respectively.

The {\em target} pendulum we wanted to control had parameters $m=1$kg, $l=0.6$m, and $dt=0.1$s. We also considered a different (i.e., {\em source}) pendulum with parameters $m = 0.5$kg, $l = $0.5m, and $dt = 0.1$s. We obtained the pmfs $\plant{k \mid k-1}$ and $\idealplant{k \mid k-1}$ from a database collected following the process from \cite{GARRABE202281}, leveraging the source code that was provided therein. We obtained $\idealcontrol{k \mid k-1}$ by controlling the source pendulum via Model Predictive Control (MPC) with a receding horizon window width of $20$ steps. The action space was $\mathcal{U}$ and the cost function at each $k$ was $\sum_{t\in k:k+H-1}(\theta_{k}^2+0.1\omega_{k}^{2})+\theta_{k+H}^{2}+0.5\omega_{k+H}^{2}$. Then, as in \cite{GARRABE202281}, we added Gaussian noise to the MPC control inputs so that $\idealcontrol{k \mid k-1}$ was $\mathcal{N}(\tilde{u}_{k},0.2)$, with $\tilde{u}_{k}$ being the  control input computed via MPC.

\begin{figure}
    \centering
    \includegraphics[width=0.8\columnwidth]{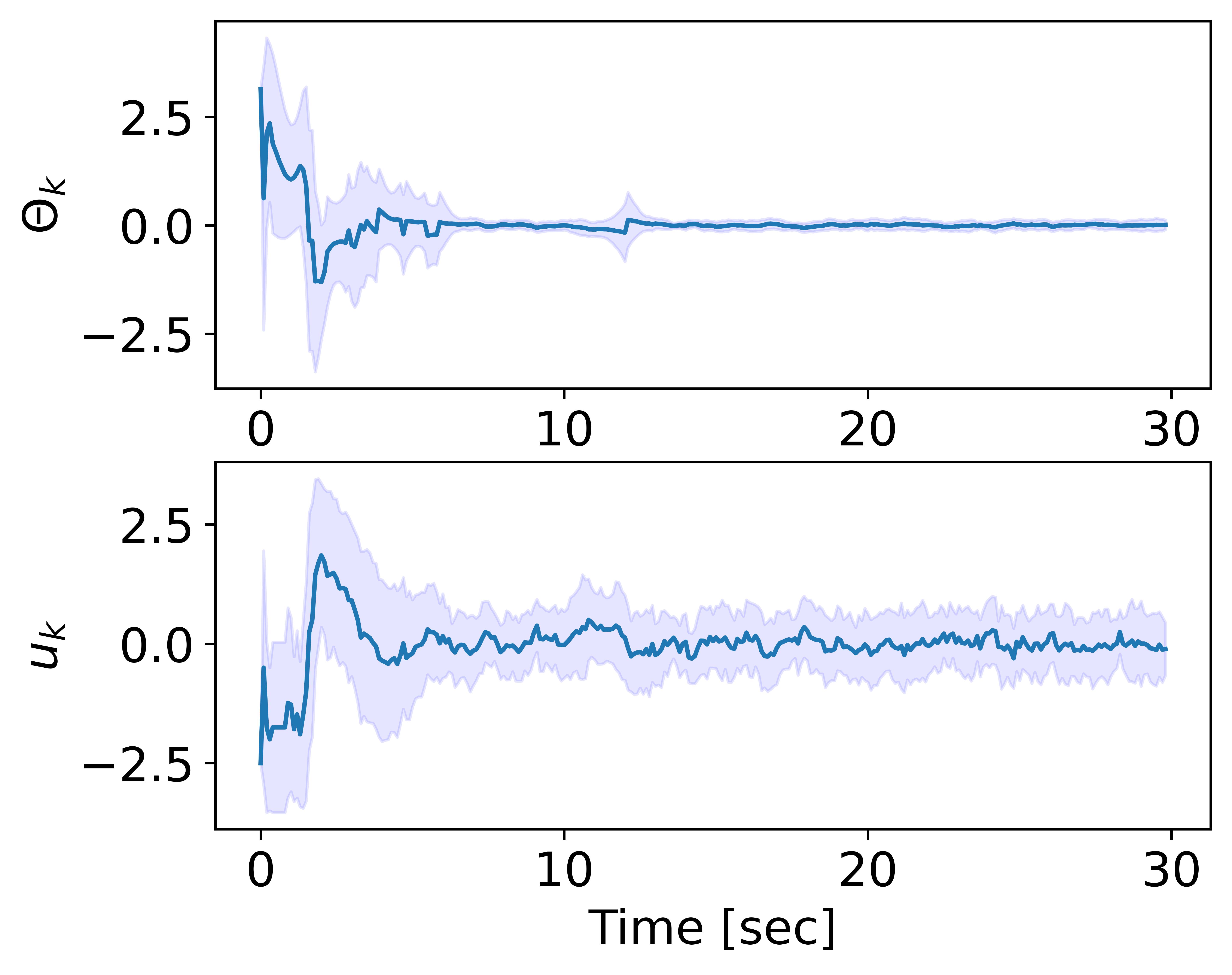}
    \caption{Target pendulum angular position and corresponding control input. Figure obtained from $20$ simulations. Bold lines represents the mean and the shaded region is confidence interval corresponding to the standard deviation.}
    \label{Simulation_1}
\end{figure}

\begin{figure}
    \centering
    \includegraphics[width=0.8\columnwidth]{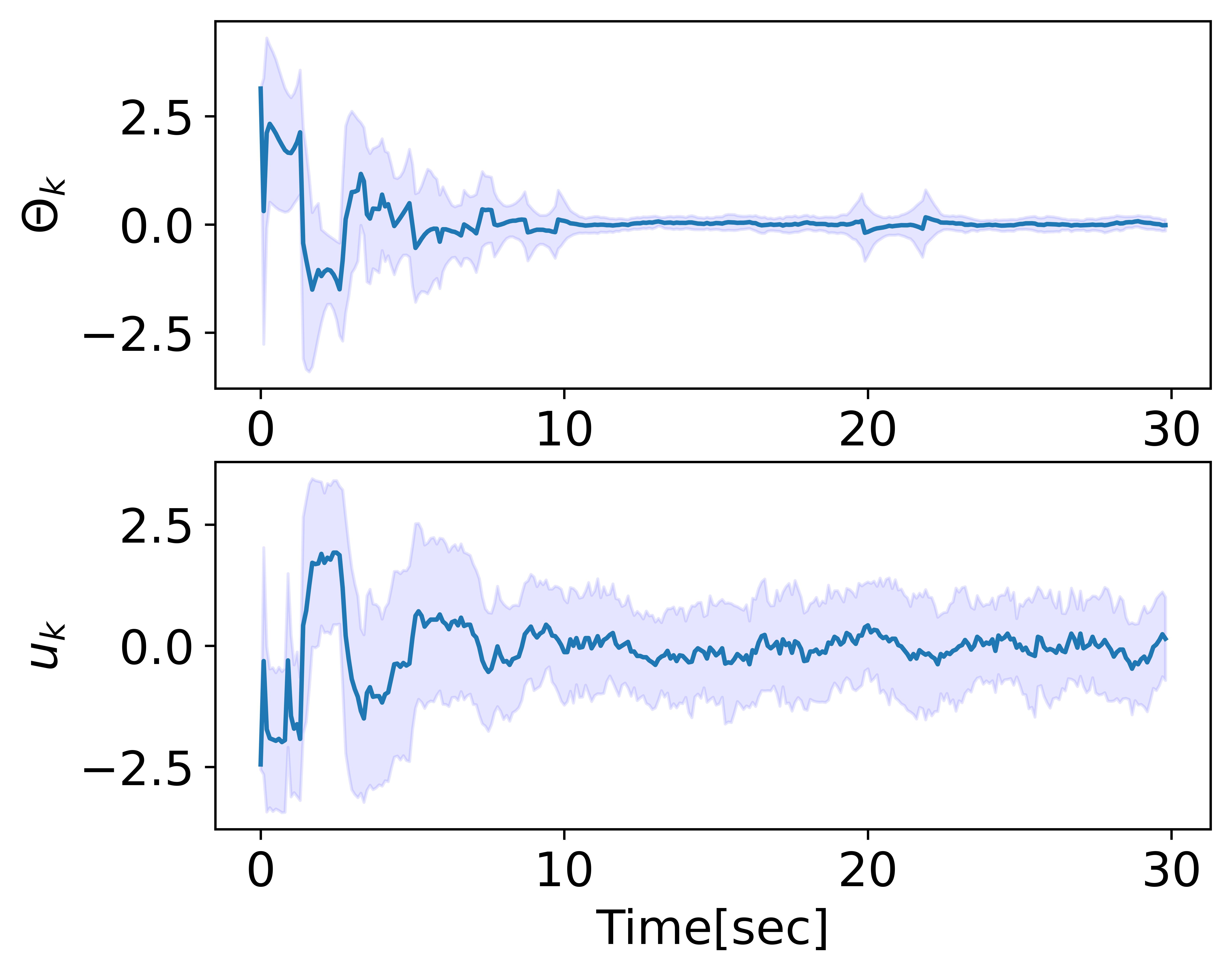}
    \caption{Angular position and control input for the target pendulum when this is controlled by a policy using the estimated cost. Figure obtained by running $20$ simulations. Bold lines and shaded area defined as in Figure \ref{Simulation_1}.}
    \label{Simulation_2}
\end{figure}

Given this set-up, we first computed  $\optimalcontrol{k \mid k-1}$ for the target pendulum using Theorem \ref{thm:prob_main} with $N=1$ and using as cost: 
\begin{equation}\label{eqn:state_cost}
    {c}(\bv{x}_{k}) = (\theta_{k} - \theta_{d})^{2}+0.01(\omega_{k} - \omega_{d})^{2},
\end{equation}
with $\theta_{d}=0$ and $\omega_{d}=0$ ($\theta_{d}=0$ corresponds to the unstable equilibrium). Then, the control input to the target pendulum was obtained by sampling from $\optimalcontrol{k \mid k-1}$. In Figure \ref{Simulation_1} the behavior is shown for the angular position of the controlled pendulum and the corresponding control input. The figure clearly illustrates that the unstable equilibrium is stabilized.  
~
Next, we illustrated the effectiveness of Algorithm \ref{alg:estimator} in reconstructing the cost given in \eqref{eqn:state_cost} that was used for policy computation. To this aim, we used a dataset of $300$ data-points collected from a single simulation where the target pendulum was controlled by the policy computed above. We defined the features as $\bv{h}(\bv{x}_{k})= \left[\mid\theta_{k}-\theta_{d}\mid,\mid\omega_{k}-\omega_{d}\mid\right]^T$. We then obtained from Algorithm \ref{alg:estimator} the weights $\bv{w}_{s}^{\ast} =[-3.3,-2.03]^{T}$ and hence the estimated cost was ${c}^{\ast}(\bv{x}_{k}) = 3.3 \mid\theta_{k}-\theta_{d}\mid + 2.03 \mid\omega_{k}-\omega_{d}\mid$. Note that the weight was higher for  $\theta_{k}$ than for $\omega_{k}$, consistently with the cost in \eqref{eqn:state_cost}. Finally, with this estimated cost, we used Theorem \ref{thm:prob_main} with $N=1$ to obtain a new policy, which we used on the target pendulum. Simulations (in Figure \ref{Simulation_2}) illustrate that this policy with the estimated cost effectively stabilizes the pendulum. 

\section{Conclusions}


We considered the problem of estimating the possibly non-convex cost of an agent by observing its interactions with a nonlinear, non-stationary, and stochastic environment. Using probabilistic descriptions from data and/or first-principles, we formulated a convex optimization problem to estimate the cost. To solve the inverse problem, we also formulated a convex finite-horizon optimal control problem and found its optimal solution. The results were turned into an algorithm and illustrated through simulations, with future work focusing on environment learning and constrained control tasks.


\bibliography{CDC_Main}
\end{document}